\documentclass[11pt, a4paper]{article}
\usepackage[cp1251]{inputenc}
\usepackage[russian]{babel}
\usepackage{amsmath} \usepackage{euscript}
\usepackage{mymatrix2}
\usepackage{longtable}
\usepackage{mathrsfs}
\usepackage{amssymb}
\begin{document}

\newcounter{bnomer} \newcounter{snomer}
\newcounter{bsnomer}
\setcounter{bnomer}{0}
\renewcommand{\thesnomer}{\thebnomer.\arabic{snomer}}
\renewcommand{\thebsnomer}{\thebnomer.\arabic{bsnomer}}
\renewcommand{\refname}{\begin{center}\large{\textbf{References}}\end{center}}

\setcounter{MaxMatrixCols}{14}

\newcommand{\sect}[1]{%
\setcounter{snomer}{0}\setcounter{bsnomer}{0}
\refstepcounter{bnomer}
\par\bigskip\begin{center}\large{\textbf{\arabic{bnomer}. {#1}}}\end{center}}
\newcommand{\sst}{%
\refstepcounter{bsnomer}
\par\bigskip\textbf{\arabic{bnomer}.\arabic{bsnomer}. }}
\newcommand{\defi}[1]{%
\refstepcounter{snomer}
\par\medskip\textbf{Definition \arabic{bnomer}.\arabic{snomer}. }{#1}\par\medskip}
\newcommand{\theo}[2]{%
\refstepcounter{snomer}
\par\textbf{Теорема \arabic{bnomer}.\arabic{snomer}. }{#2} {\emph{#1}}\hspace{\fill}$\square$\par}
\newcommand{\mtheop}[2]{%
\refstepcounter{snomer}
\par\textbf{Theorem \arabic{bnomer}.\arabic{snomer}. }{\emph{#1}}
\par\textsc{Proof}. {#2}\hspace{\fill}$\square$\par}
\newcommand{\mcorop}[2]{%
\refstepcounter{snomer}
\par\textbf{Corollary \arabic{bnomer}.\arabic{snomer}. }{\emph{#1}}
\par\textsc{Proof}. {#2}\hspace{\fill}$\square$\par}
\newcommand{\mtheo}[1]{%
\refstepcounter{snomer}
\par\medskip\textbf{Theorem \arabic{bnomer}.\arabic{snomer}. }{\emph{#1}}\par\medskip}
\newcommand{\mlemm}[1]{%
\refstepcounter{snomer}
\par\medskip\textbf{Lemma \arabic{bnomer}.\arabic{snomer}. }{\emph{#1}}\par\medskip}
\newcommand{\mprop}[1]{%
\refstepcounter{snomer}
\par\medskip\textbf{Proposition \arabic{bnomer}.\arabic{snomer}. }{\emph{#1}}\par\medskip}
\newcommand{\theobp}[2]{%
\refstepcounter{snomer}
\par\textbf{Теорема \arabic{bnomer}.\arabic{snomer}. }{#2} {\emph{#1}}\par}
\newcommand{\theop}[2]{%
\refstepcounter{snomer}
\par\textbf{Theorem \arabic{bnomer}.\arabic{snomer}. }{\emph{#1}}
\par\textsc{Proof}. {#2}\hspace{\fill}$\square$\par}
\newcommand{\theosp}[2]{%
\refstepcounter{snomer}
\par\textbf{Теорема \arabic{bnomer}.\arabic{snomer}. }{\emph{#1}}
\par\textbf{Схема доказательства}. {#2}\hspace{\fill}$\square$\par}
\newcommand{\exam}[1]{%
\refstepcounter{snomer}
\par\medskip\textbf{Example \arabic{bnomer}.\arabic{snomer}. }{#1}\par\medskip}
\newcommand{\deno}[1]{%
\refstepcounter{snomer}
\par\textbf{Definition \arabic{bnomer}.\arabic{snomer}. }{#1}\par}
\newcommand{\post}[1]{%
\refstepcounter{snomer}
\par\textbf{Предложение \arabic{bnomer}.\arabic{snomer}. }{\emph{#1}}\hspace{\fill}$\square$\par}
\newcommand{\postp}[2]{%
\refstepcounter{snomer}
\par\medskip\textbf{Proposition \arabic{bnomer}.\arabic{snomer}. }{\emph{#1}}%
\ifhmode\par\fi\textsc{Proof}. {#2}\hspace{\fill}$\square$\par\medskip}
\newcommand{\lemm}[1]{%
\refstepcounter{snomer}
\par\textbf{Lemma \arabic{bnomer}.\arabic{snomer}. }{\emph{#1}}\hspace{\fill}$\square$\par}
\newcommand{\lemmp}[2]{%
\refstepcounter{snomer}
\par\medskip\textbf{Lemma \arabic{bnomer}.\arabic{snomer}. }{\emph{#1}}
\par\textsc{Proof}. {#2}\hspace{\fill}$\square$\par\medskip}
\newcommand{\coro}[1]{%
\refstepcounter{snomer}
\par\textbf{Следствие \arabic{bnomer}.\arabic{snomer}. }{\emph{#1}}\hspace{\fill}$\square$\par}
\newcommand{\mcoro}[1]{%
\refstepcounter{snomer}
\par\textbf{Corollary \arabic{bnomer}.\arabic{snomer}. }{\emph{#1}}\par\medskip}
\newcommand{\corop}[2]{%
\refstepcounter{snomer}
\par\textbf{Следствие \arabic{bnomer}.\arabic{snomer}. }{\emph{#1}}
\par\textsc{Proof}. {#2}\hspace{\fill}$\square$\par}
\newcommand{\nota}[1]{%
\refstepcounter{snomer}
\par\medskip\textbf{Remark \arabic{bnomer}.\arabic{snomer}. }{#1}\par\medskip}
\newcommand{\propp}[2]{%
\refstepcounter{snomer}
\par\medskip\textbf{Proposition \arabic{bnomer}.\arabic{snomer}. }{\emph{#1}}
\par\textsc{Proof}. {#2}\hspace{\fill}$\square$\par\medskip}
\newcommand{\hypo}[1]{%
\refstepcounter{snomer}
\par\medskip\textbf{Conjecture \arabic{bnomer}.\arabic{snomer}. }{\emph{#1}}\par\medskip}
\newcommand{\prop}[1]{%
\refstepcounter{snomer}
\par\textbf{Proposition \arabic{bnomer}.\arabic{snomer}. }{\emph{#1}}\hspace{\fill}$\square$\par}

\newcommand{\Ind}[3]{%
\mathrm{Ind}_{#1}^{#2}{#3}}
\newcommand{\Res}[3]{%
\mathrm{Res}_{#1}^{#2}{#3}}
\newcommand{\epsi}{\epsilon}
\newcommand{\tri}{\triangleleft}
\newcommand{\Supp}[1]{%
\mathrm{Supp}(#1)}

\newcommand{\reg}{\mathrm{reg}}
\newcommand{\sreg}{\mathrm{sreg}}
\newcommand{\codim}{\mathrm{codim}\,}
\newcommand{\chara}{\mathrm{char}\,}
\newcommand{\rk}{\mathrm{rk}\,}
\newcommand{\chr}{\mathrm{ch}\,}
\newcommand{\id}{\mathrm{id}}
\newcommand{\Ad}{\mathrm{Ad}}
\newcommand{\col}{\mathrm{col}}
\newcommand{\row}{\mathrm{row}}
\newcommand{\low}{\mathrm{low}}
\newcommand{\pho}{\hphantom{\quad}\vphantom{\mid}}
\newcommand{\fho}[1]{\vphantom{\mid}\setbox0\hbox{00}\hbox to \wd0{\hss\ensuremath{#1}\hss}}
\newcommand{\wt}{\widetilde}
\newcommand{\wh}{\widehat}
\newcommand{\ad}[1]{\mathrm{ad}_{#1}}
\newcommand{\tr}{\mathrm{tr}\,}
\newcommand{\GL}{\mathrm{GL}}
\newcommand{\SL}{\mathrm{SL}}
\newcommand{\Sp}{\mathrm{Sp}}
\newcommand{\Mat}{\mathrm{Mat}}
\newcommand{\Stab}{\mathrm{Stab}}

\newcommand{\vfi}{\varphi}
\newcommand{\teta}{\vartheta}
\newcommand{\Bfi}{\Phi}
\newcommand{\Fp}{\mathbb{F}}
\newcommand{\Rp}{\mathbb{R}}
\newcommand{\Zp}{\mathbb{Z}}
\newcommand{\Cp}{\mathbb{C}}
\newcommand{\ut}{\mathfrak{u}}
\newcommand{\at}{\mathfrak{a}}
\newcommand{\nt}{\mathfrak{n}}
\newcommand{\mt}{\mathfrak{m}}
\newcommand{\htt}{\mathfrak{h}}
\newcommand{\spt}{\mathfrak{sp}}
\newcommand{\rt}{\mathfrak{r}}
\newcommand{\rad}{\mathfrak{rad}}
\newcommand{\bt}{\mathfrak{b}}
\newcommand{\gt}{\mathfrak{g}}
\newcommand{\vt}{\mathfrak{v}}
\newcommand{\pt}{\mathfrak{p}}
\newcommand{\Xt}{\mathfrak{X}}
\newcommand{\Po}{\mathcal{P}}
\newcommand{\Uo}{\EuScript{U}}
\newcommand{\Fo}{\EuScript{F}}
\newcommand{\Do}{\EuScript{D}}
\newcommand{\Eo}{\EuScript{E}}
\newcommand{\Iu}{\mathcal{I}}
\newcommand{\Mo}{\mathcal{M}}
\newcommand{\Nu}{\mathcal{N}}
\newcommand{\Ro}{\mathcal{R}}
\newcommand{\Co}{\mathcal{C}}
\newcommand{\Lo}{\mathcal{L}}
\newcommand{\Ou}{\mathcal{O}}
\newcommand{\Uu}{\mathcal{U}}
\newcommand{\Au}{\mathcal{A}}
\newcommand{\Vu}{\mathcal{V}}
\newcommand{\Bu}{\mathcal{B}}
\newcommand{\Sy}{\mathcal{Z}}
\newcommand{\Sb}{\mathcal{F}}
\newcommand{\Gr}{\mathcal{G}}
\newcommand{\rtc}[1]{C_{#1}^{\mathrm{red}}}

\author{Mikhail A. $\text{Bochkarev}^1$\and Mikhail V. $\text{Ignatyev}^2$\and Aleksandr A. $\text{Shevchenko}^3$}

\date{\small$\vphantom{1}^1$Moscow State University, Chair of higher algebra\\\texttt{mbochk@gmail.com}\\
\medskip$\vphantom{1}^2$Samara State University, Chair of algebra and geometry\\\texttt{mihail.ignatev@gmail.com}\\
\medskip$\vphantom{1}^3$Samara State University, Chair of algebra and geometry\\\texttt{shevchenko.alexander.1618@gmail.com}}
\title{\Large{Tangent cones to Schubert varieties in types $A_n$, $B_n$ and $C_n$}\mbox{$\vphantom{1}$}\footnotetext{The second and the third authors were partially supported by RFBR grants no. 13--01--97000 and 14--01--31052. The second author was partially supported by the Dynasty Foundation and by DAAD program ``Forschungsaufenthalte f$\ddot{\mathrm{u}}$r Hochschullehrer und Wissenschaftler'', ref. no. A/13/00032.}} \maketitle

\sect{Introduction and the main results}

\sst Let $G$ be a complex reductive algebraic group, $T$ a maximal torus in~$G$, $B$ a Borel subgroup in~$G$ containing $T$, and $U$ the unipotent radical of $B$. Let $\Phi$ be the root system of $G$ with respect to~$T$, $\Phi^+$ the set of positive roots with respect to $B$, $\Delta$ the set of simple roots, and $W$ the Weyl group of $\Phi$ (see \cite{Bourbaki}, \cite{Humphreys} and \cite{Humpreys2} for basic facts about algebraic groups and root systems).

Denote by $\Fo=G/B$ the flag variety and by $X_w\subseteq\Fo$ the Schubert subvariety corresponding to an element $w$ of the Weyl group $W$. Denote by $\Ou=\Ou_{p,X_w}$ the local ring at the point $p=eB\in X_w$. Let $\mt$ be the maximal ideal of~$\Ou$. The sequence of ideals $$\Ou\supseteq\mt\supseteq\mt^2\supseteq\ldots$$ is a filtration on $\Ou$. We define $R$ to be the graded algebra $$R=\mathrm{gr}\,\Ou=\bigoplus_{i\geq0}\mt^i/\mt^{i+1}.$$ By definition, the \emph{tangent cone} $C_w$ to the Schubert variety $X_w$ at the point $p$ is the spectrum of~$R$: $C_w=\mathrm{Spec}\,R$. Obviously, $C_w$ is a subscheme of the tangent space $T_pX_w\subseteq T_p\Fo$. A hard problem in studying geometry of $X_w$ is to describe $C_w$ \cite[Chapter 7]{BilleyLakshmibai}.

In 2011, D.Yu. Eliseev and A.N. Panov computed tangent cones $C_w$ for all $w\in W$ in the case $G=\mathrm{SL}_n(\mathbb{C})$, $n\leq5$ \cite{EliseevPanov}. Using their computations, A.N. Panov formulated the following Conjecture.

\hypo{\textup{(A.N. Panov, 2011)} Let $w_1$\textup{,} $w_2$ be \label{mconj}involutions\textup{,} i.e.\textup{,} $w_1^2=w_2^2=\id$. If $w_1\neq w_2$\textup{,} then $C_{w_1}\neq C_{w_2}$ as subschemes of $T_p\Fo$.}

In 2013, D.Yu. Eliseev and the second author proved that this Conjecture is true if all irreducible components of the root system $\Phi$ are of type $A_n$, $F_4$ and $G_2$ \cite{EliseevIgnatyev}. In this paper, we prove that the Conjecture is true if all irreducible components of $\Phi$ are or type $B_n$ and $C_n$. Precisely, our first main result is as follows.

\mtheo{Assume that every irreducible component of $\Phi$ is of type $B_n$ or~$C_n$\textup, $n\geq2$. Let $w_1$\textup,~$w_2$ be involutions in the Weyl group of $\Phi$ and $w_1\neq w_2$. Then the tangent cones $C_{w_1}$ and~$C_{w_2}$ do~not coincide as subschemes of $T_p\Fo$.\label{mtheo:non_red}}

Note that the similar question for the root systems $D_n$, $E_6$, $E_7$, $E_8$ remains open.

Now, let $\Au$ be the symmetric algebra of the vector space $\mt/\mt^2$, or, equivalently, the algebra of regular functions on the tangent space $T_pX_w$. Since $R$ is generated as $\Cp$-algebra by $\mt/\mt^2$, it is a~quotient ring $R=\Au/I$. By definition, the \emph{reduced tangent cone} $C_w^{\mathrm{red}}$ to $X_w$ at the point $p$ is the common zero locus in $T_pX_w$ of the polynomials $f\in I\subseteq\Au$. Clearly, if $\rtc{w_1}\neq\rtc{w_2}$, then $C_{w_1}\neq C_{w_2}$. Our second main result is as follows.

\mtheo{Assume that every irreducible component of $\Phi$ is of type $A_n$\textup, $n\geq1$\textup{,} or~$C_n$\textup, $n\geq2$. Let $w_1$\textup, $w_2$ be involutions in the Weyl group of $\Phi$ and $w_1\neq w_2$. Then the reduced tangent cones $\rtc{w_1}$ and $\rtc{w_2}$ do not coincide as subvarieties of $T_p\Fo$.\label{mtheo:red}}

Note that the similar question for other root systems remains open.

The paper is organized as follows. In the next Subsection, we introduce the main technical tool used in the proof of Theorem~\ref{mtheo:non_red}. Namely, to each element $w\in W$ one can assign a polynomial $d_w$ in the algebra of regular functions on the Lie algebra of the maximal torus $T$. These polynomials are called Kostant--Kumar polynomials \cite{KostantKumar1}, \cite{KostantKumar2}, \cite{Kumar}, \cite{Billey}. In \cite{Kumar} S. Kumar showed that if $w_1$ and~$w_2$ are arbitrary elements of~$W$ and $d_{w_1}\neq d_{w_2}$, then $C_{w_1}\neq C_{w_2}$. We give three equivalent definitions of Kostant--Kumar polynomials and formulate their properties needed for the sequel. In Subsection~\ref{sst:irred} we check that it is enough to prove Theorems~\ref{mtheo:non_red} and~\ref{mtheo:red} for irreducible root systems.

In Section~\ref{sect:non_red} we prove that if all irreducible components of $\Phi$ are of type $B_n$ and $C_n$ and $w_1$, $w_2$ are distinct involutions in $W$, then $d_{w_1}\neq d_{w_2}$, see Propositions \ref{prop:non_red_non_eq_C_1}, \ref{prop:non_red_eq_C_1}. This implies that $C_{w_1}\neq C_{w_2}$ and proves Theorem~\ref{mtheo:non_red}. The proof of Conjecture~\ref{mconj} for $A_n$, $F_4$ and $G_2$ presented in \cite{EliseevIgnatyev} is based on the similar argument.

Section~\ref{sect:red} contains the proof of Theorem~\ref{mtheo:red}. Namely, in Subsection~\ref{sst:red_coadjoint} we describe connections between the tangent cones to Schubert varieties and the geometry of coadjoint orbits of the Borel subgroup $B$. In Subsections~\ref{sst:red_A_n}, \ref{sst:red_C_n}, using the results of the second author about coadjoint orbits \cite{Ignatyev1}, \cite{Ignatyev2}, we prove Theorem~\ref{mtheo:red}, see Propositions~\ref{prop:red_A_n} and \ref{prop:red_C_n}.

\medskip\textsc{Acknowledgements}. We express our gratitude to Pro\-fes\-sor Dmitriy Timashev. His remarks helped us to bring this paper to a more clear form. Mikhail Ignatyev and Aleksandr Shevchenko were supported by RFBR grants no. 13--01--97000 and 14--01--31052; RFBR is gratefully acknowledged. Mikhail Ignatyev aknowledges support from the Dynasty Foundation. A part of this work was done during the stay of Mikhail Ignatyev at Jacobs University, Bremen, supported by DAAD program ``For\-schungsaufenthalte f$\ddot{\mathrm{u}}$r Hochschullehrer und Wissenschaftler'', ref. no. A/13/00032. Mikhail Ignatyev thanks Professor Ivan Penkov for his hospitality and useful discussions, and DAAD for the financial support.

\sst Let $w$ be an element of the Weyl group $W$. Here we give precise definition of the Kostant--Kumar polynomial $d_w$, explain how to compute it in combinatorial terms, and show that it depends only on the scheme structure of $C_w$.

The torus $T$ acts on the Schubert variety $X_w$ by left multiplications. The point $p$ is invariant under this action, hence there is the structure of a $T$-module on the local ring $\Ou$. The action of $T$ on $\Ou$ preserves the filtration by powers of the ideal $\mt$, so we obtain the structure of a $T$-module on the algebra $R=\mathrm{gr}\,\Ou$. By~\cite[Theorem~2.2]{Kumar}, $R$ can be decomposed into a direct sum of its finite-dimensional weight subspaces: $$R=\bigoplus_{\lambda\in\Xt(T)}R_{\lambda}.$$
Here $\htt$ is the Lie algebra of the torus $T$, $\Xt(T)\subseteq\htt^*$ is the character lattice of $T$ and $R_{\lambda}=\{f\in R\mid t.f=\lambda(t)f\}$ is the weight subspace of weight $\lambda$. Let $\Lambda$ be the $\Zp$-module consisting of all (possibly infinite) $\Zp$-linear combinations of linearly independent elements $e^{\lambda}$, $\lambda\in\Xt(T)$. The \emph{formal character} of~$R$ is an element of $\Lambda$ of the form $$\chr R=\sum_{\lambda\in\Xt(T)}m_{\lambda}e^{\lambda},$$ where $m_{\lambda}=\dim R_{\lambda}$.

Now, pick an element $a=\sum_{\lambda\in\Xt(T)}n_{\lambda}e^{\lambda}\in\Lambda$. Assume that there are finitely many $\lambda\in\Xt(T)$ such that $n_{\lambda}\neq0$. Given $k\geq0$, one can define the polynomial $$[a]_k=\sum_{\lambda\in\Xt(T)}n_{\lambda}\cdot\dfrac{\lambda^k}{k!}\in S=\Cp[\htt].$$ Denote $[a]=[a]_{k_0}$, where $k_0$ is minimal among all non-negative numbers $k$ such that $[a]_k\neq0$. For instance, if $a=1-e^{\lambda}$, then $[a]_0=0$ and $[a]=[a]_1=-\lambda$ (here we denote $1=e^0$).

Let $A$ be the submodule of $\Lambda$ consisting of all finite linear combinations. It is a commutative ring with respect to the multiplication $e^{\lambda}\cdot e^{\mu}=e^{\lambda+\mu}$. In fact, it is just the group ring of $\Xt(T)$. Denote the field of fractions of the ring $A$ by $Q\subseteq\Lambda$. To each element of $Q$ of the form $q=a/b$, $a$,~$b\in A$, one can assign the element $$[q]=\dfrac{[a]}{[b]}\in\Cp(\htt)$$ of the field of rational functions on $\htt$. Note that this element is well-defined~\cite{Kumar}.

There exists an involution $q\mapsto q^*$ on $Q$ defined by $$e^{\lambda}\mapsto (e^{\lambda})^*=e^{-\lambda}.$$ It turns out \cite[Theorem 2.2]{Kumar} that the character $\chr R$ belongs to $Q$, hence $(\chr R)^*\in Q$, too. Finally, we put $$c_w=[(\chr R)^*],\ d_w=(-1)^{l(w)}\cdot c_w\cdot\prod_{\alpha\in\Phi^+}\alpha.$$ Here $l(w)$ is the length of $w$ in the Weyl group $W$ with respect to the set of simple $\Delta$. Evidently, $c_w$ and $d_w$ belong to $\Cp(\htt)$; in fact, $d_w$ is a polynomial, i.e., it belongs to the algebra $S=\Cp[\htt]$ of regular functions on $\htt$, see \cite{KostantKumar2} and \cite[Theorem 7.2.6]{BilleyLakshmibai}. \defi{Let $w$ be an element of the Weyl group $W$. The polynomial $d_w\in S$ is called the \emph{Kostant--Kumar polynomial} associated with $w$.}

It follows from the definition that $c_w$ and $d_w$ depend only on the canonical structure of a $T$-module on the algebra $R$ of regular functions on the tangent cone $C_w$. Thus, to prove that the tangent cones corresponding to elements $w$, $w'$ of the Weyl group are distinct, it is enough to check that ${c_{w_1}\neq c_{w_2}}$, or, equivalently, $d_{w_1}\neq d_{w_2}$.

On the other hand, there is a purely combinatorial description of Kostant--Kumar polynomials. To give this description, we need some more notation. Let $w$, $v$ be elements of~$W$. Fix a reduced decomposition of the element $w=s_{i_1}\ldots s_{i_l}$. (Here $\alpha_1,\ldots,\alpha_n\in\Delta$ are simple roots and $s_i$ is the simple reflection corresponding to $\alpha_i$.) Put
\begin{equation*}
c_{w,v}=(-1)^{l(w)}\cdot\sum\dfrac{1}{s_{i_1}^{\epsi_1}\alpha_{i_1}}\cdot\dfrac{1}{s_{i_1}^{\epsi_1}
s_{i_2}^{\epsi_2}\alpha_{i_2}}\cdot\ldots\cdot\dfrac{1}{s_{i_1}^{\epsi_1}\ldots s_{i_l}^{\epsi_l}\alpha_{i_l}},
\end{equation*}
where the sum is taken over all sequences $(\epsi_1,\ldots,\epsi_l)$ of zeroes and units such that $s_{i_1}^{\epsi_1}\ldots s_{i_l}^{\epsi_l}=v$. Actually, the element $c_{w,v}\in\Cp(\htt)$ depends only on $w$ and $v$, not on the choice of a reduced decomposition of $w$ \cite[Section 3]{Kumar}.

\exam{Let $\Phi=A_n$. Put $w=s_1s_2s_1$. To compute $c_{w,\id}$, we should take the sum over two sequences, $(0,0,0)$ and $(1,0,1)$. Hence
\begin{equation*}
c_{w,\id}=(-1)^3\cdot\left(\dfrac{1}{\alpha_1\alpha_2\alpha_1}+\dfrac{1}{-\alpha_1(\alpha_1+\alpha_2)\alpha_1}\right)
=\dfrac{1}{\alpha_1\alpha_2(\alpha_1+\alpha_2)}.
\end{equation*}}

A remarkable fact is that $c_{w,\id}=c_w$, hence to prove that the tangent cones to Schubert varieties do not coincide as subschemes, we need only combinatorics of the Weyl group. Note also that for classical Weyl groups, elements $c_{w,v}$ are closely related to Schubert polynomials \cite{Billey}.

Finally, we will present an original definition of elements $c_{w,v}$ using so-called nil-Hecke ring (see \cite{Kumar} and~\cite[Section 7.1]{BilleyLakshmibai}). The group $W$ naturally acts on $\Cp(\htt)$ by automorphisms. Denote by $Q_W$ the vector space over $\Cp(\htt)$ with basis $\{\delta_w,\ w\in W\}$. It is a ring with respect to the multiplication $$f\delta_v\cdot g\delta_w=fv(g)\delta_{vw}.$$
This ring is called the \emph{nil-Hecke ring}. To each $i$ from 1 to $n$ put $$x_i=\alpha_i^{-1}(\delta_{s_i}-\delta_{\id}).$$
Let $w\in W$ and $w=s_{i_1}\ldots s_{i_l}$ be a reduced decomposition of $W$. Then the element $$x_w=x_{i_1}\ldots x_{i_l}$$ does not depend on the choice of a reduced decomposition \cite[Proposition 2.1]{KostantKumar1}.

Moreover, it turns out that $\{x_w,\ w\in W\}$ is a $\Cp(\htt)$-basis of~$Q_W$ \cite[Proposition 2.2]{KostantKumar1}, and
\begin{equation*}
x_w=\sum\nolimits_{v\in W}c_{w,v}\delta_v.\\
\end{equation*}
Actually, if $w,v\in W$, then
\begin{equation}
\begin{split}
&\text{a) }x_v\cdot x_w=\begin{cases}x_{vw},&\text{ if }l(vw)=l(v)+l(w),\\
0,&\text{ otherwise},
\end{cases}\\
&\text{b) }c_{w,v}=-v(\alpha_i)^{-1}(c_{ws_i,v}+c_{ws_i,vs_i}),\text{ if }l(ws_i)=l(w)-1,\\
&\text{c) }c_{w,v}=\alpha_i^{-1}(s_i(c_{s_iw,s_iv})-c_{s_iw,v}),\text{ if }l(s_iw)=l(w)-1.\\
\end{split}\label{formula_x_v_x_w}
\end{equation}
The first property is proved in \cite[Proposition 2.2]{KostantKumar1}. The second and the third properties follow immediately from the first one and the definitions (see also the proof of \cite[Corollary 3.2]{Kumar}).

\sst \label{sst:irred}In this Subsection we check that it is enough to prove Theorems~\ref{mtheo:non_red} and~\ref{mtheo:red} for irreducible root systems. Suppose~$\Phi$ is a union of its subsystems $\Phi_1$ and $\Phi_2$ contained in mutually orthogonal subspaces. Let $W_1$, $W_2$ be the Weyl groups of $\Phi_1$, $\Phi_2$ respectively, so $W=W_1\times W_2$. Denote $\Delta_1=\Delta\cap\Phi_1=\{\alpha_1,\ldots,\alpha_r\}$ and $\Delta_2=\Delta\cap\Phi_2=\{\beta_1,\ldots,\beta_s\}$, then $$\Cp[\htt]\cong\Cp[\alpha_1,\ldots,\alpha_r,\beta_1,\ldots,\beta_s].$$

Given $v\in W_1$, denote by $d_v^1$ its Kostant--Kumar polynomial. We can consider $d_v^1$ as an element of~$\Cp(\htt)$ depending only on $\alpha_1,\ldots,\alpha_r$. We define $c_v^1\in\Cp(\htt)$ by the similar way. Given $v\in W_2$, we define $d_v^2\in\Cp[\htt]$ and $c_v^2\in\Cp(\htt)$; they depend only on $\beta_1,\ldots,\beta_s$.

\propp{Let $w\in W$, $w_1\in W_1$, $w_2\in W_2$ and $w=w_1w_2$. \label{prop:irred_non_red}Then $$d_w=d_{w_1}^1d_{w_2}^2,\ c_w=c_{w_1}^1c_{w_2}^2.$$}{One can repeat literally the proof of \cite[Proposition 1.6]{EliseevIgnatyev} to obtain the result. Namely, by~$s_i$ (resp. by~$r_j$), we denote the simple reflection corresponding to a simple root $\alpha_i$ (resp. $\beta_j$). Let $l_i$ be the length function on $W_i$ with respect to $\Delta_i$, $i=1,2$. It is well-known that $l_i(v)=l(v)$ for all $v\in W_i$. Hence if
\begin{equation*}
w_1=s_{i_1}\ldots s_{i_p},\ w_2=r_{j_1}\ldots r_{j_q}
\end{equation*}
are reduced decompositions for $w_i$ in $W_i$, then they are reduced decompositions for $w_i$ in $W$. Moreover, $$l(w)=l(w_1)+l(w_2)=l_1(w_1)+l_2(w_2).$$ This means that
$$w=s_{i_1}\ldots s_{i_p}r_{j_1}\ldots r_{j_q}$$
is a reduced decomposition for $w$ in $W$.

It follows from $W=W_1\times W_2$ that $$s_{i_1}^{\epsi_1}\ldots s_{i_p}^{\epsi_p}r_{j_1}^{\delta_1}\ldots r_{j_q}^{\delta_q}=\id,$$ $\epsi_i,\delta_j\in\{0,1\}$, is equivalent to $$s_{i_1}^{\epsi_1}\ldots s_{i_p}^{\epsi_p}=r_{j_1}^{\delta_1}\ldots r_{j_q}^{\delta_q}=\id.$$ Since all $s_i$'s (resp. $r_j$'s) act identically on $\Phi_2$ (resp. on $\Phi_1$), we obtain
\begin{equation*}
\begin{split}
c_w=(-1)^{l_1(w_1)+l_2(w_2)}&\cdot\sum\left(\dfrac{1}{s_{i_1}^{\epsi_1}\alpha_{i_1}}\cdot\dfrac{1}{s_{i_1}^{\epsi_1}
s_{i_2}^{\epsi_2}\alpha_{i_2}}\cdot\ldots\cdot\dfrac{1}{s_{i_1}^{\epsi_1}\ldots s_{i_p}^{\epsi_p}\alpha_{i_p}}\vphantom{\dfrac{1}{r_{j_1}^{\delta_1}\ldots r_{j_q}^{\delta_q}\beta_{j_q}}}\right.\\
&\times\left.\dfrac{1}{r_{j_1}^{\delta_1}\beta_{j_1}}\cdot\dfrac{1}{r_{j_1}^{\delta_1}
r_{j_2}^{\delta_2}\beta_{j_2}}\cdot\ldots\cdot\dfrac{1}{r_{j_1}^{\delta_1}\ldots r_{j_q}^{\delta_q}\beta_{j_q}}\right)=c_{w_1}^1c_{w_2}^2.
\end{split}
\end{equation*}
The second equality is proved. The first equality follows immediately from the second one and the obvious fact that $\Phi^+=\Phi_1^+\cup\Phi_2^+$.}

Thus, to prove Theorem~\ref{mtheo:non_red}, it suffice to check it for irreducible root systems of types $B_n$ and $C_n$, because $\Cp[\htt]$ is a unique factorization domain.

Now, let $G\cong G_1\times G_2$, where $G_1$, $G_2$ are reductive subgroups of $G$, $T_i=T\cap G_i$ is a maximal torus in $G_i$, $i=1,2$, and the root system of $G_i$ with respect to $T_i$ is isomorphic to $\Phi_i$. Then $B_i=B\cap G_i$ is a Borel subgroup in $G_i$ containing $T_i$. Denote by $\Fo_i=G_i/B_i$ the corresponding flag variety. Then $$\Fo=\Fo_1\times\Fo_2$$ and $$T_p\Fo=T_p\Fo_1\times T_p\Fo_2$$ as algebraic varieties. If $w\in W$ and $w=w_1w_2$, $w_i\in W_i$, $i=1,2$, then $$\rtc{w}\cong\rtc{w_1,G_1}\times\rtc{w_2,G_2}$$ as affine varieties. Here $\rtc{w_i,G_i}$, $i=1,2$, denotes the tangent cone to the Schubert subvariety $X_{w_i}$ of the flag variety $\Fo_i$. Furthermore, note that $w$ is an involution if and only if $w_1$ and $w_2$ are involutions, too. This means that it suffice to prove that Theorem~\ref{mtheo:red} holds for all irreducible root systems of types $A_n$ and $C_n$.

\sect{Types $B_n$ and $C_n$, non-reduced case}\label{sect:non_red}

\sst Throughout this Section, $\Phi$ denotes an irreducible root system of type $B_n$ or $C_n$. In this Subsection, we briefly recall some basic facts about $\Phi$. Let $\epsi_1$, $\ldots$, $\epsi_n$ be the standard basis of the Euclidean space~$\Rp^n$. As usual, we identify the set $\Phi^+$ of positive roots with one of the following subsets of $\Rp^n$:
\begin{equation*}
\begin{split}
B_n^+&=\{\epsi_i-\epsi_j,~\epsi_i+\epsi_j,~1\leq i<j\leq n\}\cup\{\epsi_i,~1\leq i\leq n\},\\
C_n^+&=\{\epsi_i-\epsi_j,~\epsi_i+\epsi_j,~1\leq i<j\leq n\}\cup\{2\epsi_i,~1\leq i\leq n\},
\end{split}
\end{equation*}
so $W$ can be considered as a subgroup of the orthogonal group $O(\Rp^n)$.

Let $S_{\pm n}$ denote the symmetric group on $2n$ letters $1,\ldots,n,-n,\ldots,-1$. The Weyl group $W$ is isomorphic to the \emph{hyperoctahedral group}, that is, the subgroup of $S_{\pm n}$ consisting of permutations ${w\in S_{\pm n}}$ such that $w(-i)=-w(i)$ for all $1\leq i\leq n$. The isomorphism is given by
\begin{equation*}
\begin{split}
&s_{\epsi_i-\epsi_j}\mapsto(i,j)(-i,-j),\\
&s_{\epsi_i+\epsi_j}\mapsto(i,-j)(-i,j),\\
&s_{\epsi_i}=s_{2\epsi_i}\mapsto(i,-i).
\end{split}
\end{equation*}
In the sequel, we will identify $W$ with the hyperoctahedral group.

\nota{i) Note that every $w\in W$ is completely determined by its restriction to  the subset $\{1,\ldots,n\}$. This allows us to use the usual two-line notation: if $w(i)=w_i$ for $1\leq i\leq n$, then we will write $w=\begin{pmatrix}1&2&\ldots&n\\w_1&w_2&\ldots&w_n\end{pmatrix}$. For instance, if $\Phi=C_5$, then $$s_{\epsi_1+\epsi_5}s_{2\epsi_3}s_{\epsi_2-\epsi_4}=\begin{pmatrix}1&2&3&4&5\\-5&4&-3&2&-1\end{pmatrix}.$$

ii) Note also that the set of simple roots has the following form: $\Delta=\{\alpha_1,\ldots,\alpha_n\}$, where $\alpha_1=\epsi_1-\epsi_2$, $\ldots$, $\alpha_{n-1}=\epsi_{n-1}-\epsi_n$, and
\begin{equation*}
\alpha_n=\begin{cases}\epsi_n,&\text{if }\Phi=B_n,\\
2\epsi_n,&\text{if }\Phi=C_n.
\end{cases}
\end{equation*}
We will denote the corresponding simple reflections by $s_1,\ldots,s_n$.}

We say that $v$ is less or equal to $w$ with respect to the \emph{Bruhat order}, written $v\leq w$, if some reduced decomposition for $v$ is a subword of some reduced decomposition for $w$. It is well-known that this order plays the crucial role in many geometric aspects of theory of algebraic groups. For instance, the Bruhat order encodes the incidences among Schubert varieties, i.e., $X_v$ is contained in $X_w$ if and only if $v\leq w$. It turns out that $c_{w,v}$ is non-zero if and only if $v\leq w$ \cite[Corollary 3.2]{Kumar}. For example, $c_w=c_{w,\id}$ is non-zero for \emph{any} $w$, because $\id$ is the smallest element of $W$ with respect to the Bruhat order. Note that given $v,w\in W$, there exists $g_{w,v}\in S=\Cp[\htt]$ such that
\begin{equation}
c_{w,v}=g_{w,v}\cdot\prod_{\alpha>0,~s_{\alpha}v\leq w}\alpha^{-1},\label{formula:dyer}
\end{equation}
see \cite{Dyer} and \cite[Theorem 7.1.11]{BilleyLakshmibai}

There exists a nice combinatorial description of the Bruhat order on the hyperoctahedral group. Namely, given $w\in W$, denote by $X_w$ the $2n\times 2n$ matrix of the form
\begin{equation*}
(X_w)_{i,j}=\begin{cases}1,&\text{ if }w(j)=i,\\
0&\text{ otherwise}.
\end{cases}
\end{equation*}
The rows and the columns of this matrix are indicated by the numbers $1,\ldots,n,-n,\ldots,1$. It is called the $0$--1 matrix, permutation matrix or rook placement for $w$. Define the matrix $R_w$ by putting its $(i,j)$th element to be equal to the rank of the lower left $(n-i+1)\times j$ submatrix of $X_w$. In other words, $(R_w)_{i,j}$ is just the number or rooks located non-strictly to the South-West from $(i,j)$.

\exam{Let $n=4$, $w=\begin{pmatrix}1&2&3&4\\-3&-2&4&-1\end{pmatrix}$. Here we draw the matrices $X_w$ and $R_w$ (rooks are marked by $\otimes$):
\begin{equation*}X_w=
\mymatrix{
\pho& \pho& \pho& \pho& \otimes& \pho& \pho& \pho\\
\pho& \pho& \pho& \pho& \pho& \pho& \otimes& \pho\\
\pho& \pho& \pho& \pho& \pho& \pho& \pho& \otimes\\
\pho& \pho& \otimes& \pho& \pho& \pho& \pho& \pho\\
\pho& \pho& \pho& \pho& \pho& \otimes& \pho& \pho\\
\otimes& \pho& \pho& \pho& \pho& \pho& \pho& \pho\\
\pho& \otimes& \pho& \pho& \pho& \pho& \pho& \pho\\
\pho& \pho& \pho& \otimes& \pho& \pho& \pho& \pho\\}
\ ,\ R_w=
\mymatrix{
\fho1& \fho2& \fho3& \fho4& \fho5& \fho6& \fho7& \fho8\\
\fho1& \fho2& \fho3& \fho4& \fho4& \fho5& \fho6& \fho7\\
\fho1& \fho2& \fho3& \fho4& \fho4& \fho5& \fho5& \fho6\\
\fho1& \fho2& \fho3& \fho4& \fho4& \fho5& \fho5& \fho5\\
\fho1& \fho2& \fho2& \fho3& \fho3& \fho4& \fho4& \fho4\\
\fho1& \fho2& \fho2& \fho3& \fho3& \fho3& \fho3& \fho3\\
\fho0& \fho1& \fho1& \fho2& \fho2& \fho2& \fho2& \fho2\\
\fho0& \fho0& \fho0& \fho1& \fho1& \fho1& \fho1& \fho1\\
}\ .
\end{equation*}\label{exam:Bruhat}}

Let $X$ and $Y$ be matrices with integer entries. We say that $X\leq Y$ if $X_{i,j}\leq Y_{i,j}$ for all $i,j$. It turns out that given $v$, $w\in W$,
\begin{equation}
v\leq w \text{ if and only if }R_v\leq R_w\label{formula:Bruhat}
\end{equation}(see, e.g., \cite{Proctor} or \cite[Theorem 8.1.8]{BjornerBrenti}).

\sst\label{sst:support} In this Subsection, we introduce some more notation and prove technical, but crucial\break Lemma~\ref{lemm:crucial_non_red}. We define the maps $\row\colon\Phi^+\to\Zp$ and $\col\colon\Phi^+\to\Zp$ by
\begin{equation*}
\begin{split}
&\row(\epsi_i-\epsi_j)=j,~\row(\epsi_i+\epsi_j)=-j,~\row(\epsi_i)=0,~\row(2\epsi_i)=-i,\\
&\col(\epsi_i-\epsi_j)=\col(\epsi_i+\epsi_j)=\col(\epsi_i)=\col(2\epsi_i)=i.
\end{split}
\end{equation*}
For any $k$ from $-n$ to $n$, put
\begin{equation*}
\begin{split}
\Ro_k&=\{\alpha\in\Phi^+\mid\row(\alpha)=k\},\\
\Co_k&=\{\alpha\in\Phi^+\mid\col(\alpha)=k\}.\\
\end{split}
\end{equation*}
The set $\Ro_k$ (resp. $\Co_k$) is called the $k$th \emph{row} (resp. the $k$th \emph{column}) of $\Phi^+$.

\defi{Let $\sigma\in W$ be an involution.\label{defi:support} We define the \emph{support} $\Supp{\sigma}$ of the involution $\sigma$ by the following rule:
\begin{equation*}
\begin{split}
&\text{if }1\leq i<j\leq n\text{ and }\sigma(i)=j,\text{ then }\epsi_i-\epsi_j\in\Supp{\sigma},\\
&\text{if }1\leq i<j\leq n\text{ and }\sigma(i)=-j,\text{ then }\epsi_i+\epsi_j\in\Supp{\sigma},\\
&\text{if }\Phi=B_n,~1\leq i\leq n\text{ and }\sigma(i)=-i,\text{ then }\epsi_i\in\Supp{\sigma},\\
&\text{if }\Phi=C_n,~1\leq i\leq n\text{ and }\sigma(i)=-i,\text{ then }2\epsi_i\in\Supp{\sigma}.\\
\end{split}
\end{equation*}
By definition, $\Supp{\sigma}$ is an orthogonal subset of $\Phi^+$. Note that $$\sigma=\prod_{\beta\in\Supp{\sigma}}s_{\beta},$$ where the product is taken in any fixed order. Note also that for any $k$ one has $$|\Supp{\sigma}\cap\Co_k|\leq1.$$}

\exam{Let $\Phi=C_6$ and $\sigma=\begin{pmatrix}1&2&3&4&5&6\\-6&-2&5&4&3&-1\end{pmatrix}$. Then $$\Supp{\sigma}=\{\epsi_1+\epsi_6,2\epsi_2,\epsi_3-\epsi_5\}.$$}

\nota{i) Denote the set of involutions by $I(W)$. By \cite{Ignatyev2}, if $\sigma,\tau\in I(W)$, then
\begin{equation}
\sigma\leq\tau\text{ if and only if }R_{\sigma}^*\leq R_{\tau}^*,\label{formula:Bruhat_invol}
\end{equation}
where $R_w^*$ is the strictly lower-triangular part of $R_w$, i.e.,
\begin{equation*}
(R_w^*)_{i,j}=\begin{cases}
(R_w)_{i,j}&\text{if }i>j,\\
0,&\text{if }i\leq j.
\end{cases}
\end{equation*}

ii) Using Formulas (\ref{formula:Bruhat}) or (\ref{formula:Bruhat_invol}), one can easily check that if $\alpha\in\Co_1$ and $\beta\notin\Co_1$, then $s_{\alpha}\nleq s_{\beta}$. One can also check that
\begin{equation*}
s_{\epsi_1-\epsi_2}<\ldots<s_{\epsi_1-\epsi_n}<s_{\epsi_1+\epsi_n}<\ldots<s_{\epsi_1+\epsi_2}.
\end{equation*}
Further, $s_{\epsi_i-\epsi_j}<s_{2\epsi_k}$, but $s_{\epsi_i+\epsi_j}\nless s_{2\epsi_k}$, $s_{2\epsi_k}\nless s_{\epsi_i+\epsi_j}$ for all $i,j,k$.}

The following Lemma plays the crucial role in the proof of Theorem~\ref{mtheo:non_red} (cf. \cite[Lemmas 2.4, 2.5]{EliseevIgnatyev}).

\lemmp{Let $w\in W$ be an involution. If $\Supp{w}\cap\Co_1=\varnothing$\textup{,} then $\alpha$ divides $d_w$ in the polynomial ring $\Cp[\htt]$ for all $\alpha\in\Co_1$. If $\Supp{w}\cap\Co_1=\{\beta\}$\textup{,} then $\beta$ does not divide $d_w$ in $\Cp[\htt]$.\label{lemm:crucial_non_red}}{Denote by $\wt W$ the subgroup of $W$ generated by $s_2,\ldots,s_n$. Suppose $\Supp{w}\cap\Co_1=\varnothing$, then $w\in\wt W$. Denote by $\wt\Phi$ the root system corresponding to $\wt W$; in fact, $\wt\Phi^+=\Phi^+\setminus\Co_1$.

Let $\wt d_w\in\wt S=\Cp[\alpha_2,\ldots,\alpha_n]$ be the Kostant--Kumar polynomial of $w$ considered as an element of~$\wt W$; define $\wt c_w\in\Cp(\alpha_2,\ldots,\alpha_n)$ by the similar way. Since $\wt W$ is a parabolic subgroup of $W$, the length of $w$ as an element of $\wt W$ equals the length of $w$ as an element of $W$. Further, any reduced decomposition for $w$ in $\wt W$ is a reduced decomposition for $w$ in $W$. This means that $\wt c_w=c_w$, so
\begin{equation*}
\begin{split}
d_w=(-1)^{l(w)}\cdot\prod_{\alpha\in\Phi^+}\alpha\cdot c_w=(-1)^{l(w)}\cdot\prod_{\alpha\in\Co_1}\alpha\cdot\prod_{\alpha\in\wt\Phi^+}\alpha\cdot\wt c_w=\wt d_w\cdot\prod_{\alpha\in\Co_1}\alpha.
\end{split}
\end{equation*}
In particular, $\alpha$ divides $d_w$ for all $\alpha\in\Co_1$.

Now, suppose $\Supp{w}\cap\Co_1=\{\beta\}$. By \cite[Proposition 1.10]{Humpreys2}, there exists a unique $v\in\wt W$ such that $w=uv$ and $l(us_i)=l(u)+1$ for all $2\leq i\leq n$ (or, equivalently, $u(\alpha_i)>0$ for all $2\leq i\leq n$). Furthermore, $l(w)=l(u)+l(v)$. One can easily check that
\begin{equation*}
\begin{split}
&\text{if }\beta=\epsi_1-\epsi_j\text{ (i.e., $w(1)=j$), then}\\
&u=\begin{pmatrix}1&2&3&\ldots&j-1&j&j+1&\ldots\\j&1&2&\ldots&j-2&j-1&j+1&\ldots\end{pmatrix}=s_{j-1}\ldots s_2s_1,\\
&\text{if }\beta=\epsi_1+\epsi_j\text{ (i.e., $w(1)=-j$), then}\\
&u=\begin{pmatrix}1&2&3&\ldots&j-1&j&j+1&\ldots\\-j&1&2&\ldots&j-2&j-1&j+1&\ldots\end{pmatrix}=s_js_{j+1}\ldots s_{n-1}s_n s_{n-1}\ldots s_2s_1,\\
&\text{if }\Phi=B_n\text{ and }\beta=\epsi_1,\text{ or }\Phi=C_n\text{ and }\beta=2\epsi_1\text{ (i.e., $w(1)=-1$), then}\\
&u=\begin{pmatrix}1&2&3&\ldots\\-1&2&3&\ldots\end{pmatrix}=s_1s_2\ldots s_{n-1}s_n s_{n-1}\ldots s_2s_1.
\end{split}
\end{equation*}
Let us consider these three cases separately.

i) Suppose $\beta=\epsi_1-\epsi_j$. Note that $W$ acts on $\Cp(\htt)$ by automorphisms. Using (\ref{formula_x_v_x_w}) and arguing as in the proof of
\cite[Lemma 2.5]{EliseevIgnatyev}, one can easily show that
\begin{equation*}
c_w=-\dfrac{c_{us_1,g_0}g_0(c_{v,g_0^{-1}})}{\beta}-\sum_{g\leq u,\ g^{-1}\leq v,\ g\neq g_0}\dfrac{c_{us_1,g}g(c_{v,g^{-1}})}{g(\alpha_1)}=
\beta^{-1}\cdot g_0(c_{v,g_0^{-1}})\cdot\dfrac{K}{L}+\dfrac{M}{N}
\end{equation*}
(cf. Formula (7) from \cite{EliseevIgnatyev}). Here
$$g_0=us_1=\begin{pmatrix}1&2&3&\ldots&j-1&j&j+1&\ldots\\1&j&2&\ldots&j-2&j-1&j+1&\ldots\end{pmatrix}=s_{j-1}\ldots s_3s_2,$$
and $K,L,M,N\in\Cp[\htt]$ are coprime polynomials such that $\beta$ divides neither $K$ nor $N$. We see that it is enough to check that $v^{-1}\geq g_0$, or, equivalently, $R_{v^{-1}}\geq R_{g_0}$. For $j=2$, there is nothing to prove, so we may assume $j>2$. Note that
\begin{equation*}
(R_{g_0})_{p,q}=\begin{cases}
1,&\text{if }2\leq q<p\leq j,\\
q-p+2,&\text{if }2\leq-p<-q\leq j,\\
(R_{\id})_{p,q}&\text{otherwise}.
\end{cases}
\end{equation*}
Since $\id$ is the smallest element of $W$ with respect to the Bruhat order, it remains to show that $(R_{v^{-1}})_{p,q}\geq(R_{g_0})_{p,q}$, if $2\leq q<p\leq j$ or $2\leq-p<-q\leq j$. But $v^{-1}(2)=wu(2)=w(1)=j$, hence if $2\leq q<p\leq j$, then $(R_{v^{-1}})_{p,q}\geq1$, Thus, it suffices to check that
\begin{equation}
(R_{v^{-1}})_{p,q}\geq(R_{g_0})_{p,q},\text{ if }2\leq-p<-q\leq j.\label{formula:last_cond_i_m_j}
\end{equation}

Let us write $1\prec 2\prec\ldots\prec n\prec -n\prec\ldots\prec-2\prec-1$. Denote $[\pm n]=\{1,2,\ldots,n,-n,\ldots,-2,-1\}$. For $a,b\in[\pm n]$, we will write $a\preceq b$ if $a\prec b$ or $a=b$. Further, for $(a,b),(c,d)\in[\pm n]\times[\pm n]$, we will write $(a,b)\preceq(c,d)$ if $c\preceq a$ and $b\preceq d$. Given $(a,b)\in[\pm n]\times[\pm n]$, put
$$N_{a,b}=\{(c,d)\in[\pm n]\times[\pm n]\mid(c,d)\preceq(a,b)\}.$$
Define the sets $\Au,\Bu\subset\Zp\times\Zp$ as follows:
\begin{equation*}
\begin{split}
\Au&=\bigcup_{3\leq a\leq j}\{(-(a-1),-a)\}=\bigcup_{3\leq a\leq j}\{(g_0(-a),-a)\},\\
\Bu&=\bigcup_{b\in[\pm n]}\{(v^{-1}(b), b)\}.
\end{split}
\end{equation*}
To check (\ref{formula:last_cond_i_m_j}), it is enough to construct an injective map $\phi\colon\Au\to\Bu$ such that
\begin{equation}
|\Au\cap N_{p,q}|\leq|\phi(\Au)\cap N_{p,q}|\text{ for all }2\leq-p<-q\leq j.\label{formula:Au_N}
\end{equation}

To do this, we need some more notation. Put
\begin{equation*}
\begin{split}
&Z=\{3,4,\ldots,j\},~Z'=u(Z)=\{2,3,\ldots,j-1\},\\
&Z_1=\{k\in Z'\mid w(k)=k\},~Z_2'=\{k\in Z'\mid w(k)=-k\},\\
&Z_3'=\{k\in Z'\mid w(k)=l>j\},~Z_4'=\{k\in Z'\mid w(k)=-l,l>j\},\\
&Z_5'=\{k,l\in Z'\mid k<l,w(k)=l\},~Z_6'=\{k,l\in Z'\mid k<l,w(k)=-l\},\\
&Z_i=u^{-1}(Z_i'),~\Au_i=\bigcup_{a\in Z_i}\{(-(a-1),-a)\},~1\leq i\leq 6.
\end{split}
\end{equation*}
Note that $\Au$ is a disjoint union of $\Au_1,\ldots,\Au_6$, because $Z'$ is a disjoint union of $Z_1',\ldots,Z_6'$. If $a\in Z_1$, then $v^{-1}(-a)=wu(-a)=w(-(a-1))=-(a-1)$. In this case, we put $\phi(-(a-1),-a)=(-(a-1),-a)$. If $a\in Z_2$, then $v^{-1}(a)=w(a-1)=-(a-1)$. Here we put $\phi(-(a-1),-a)=(-(a-1),a)$.

If $a\in Z_3$, then $v^{-1}(a)=w(a-1)=l$ for some $l>j$. Hence $v^{-1}(-l)=w(-l)=-(a-1)$. (Note that $-l\prec-a$.) In this case, we put $\phi(-(a-1),a)=(-(a-1),-l)$. If $a\in Z_4$, then $v^{-1}(a)=w(a-1)=-l$ for some $l>j$. Hence $v^{-1}(l)=w(l)=-(a-1)$. Here we put $\phi(-(a-1),a)=(-(a-1),l)$.

If $k,l\in Z_5$, $k<l$ and $w(k-1)=l-1$, then, clearly, $v^{-1}(-k)=w(-(k-1))=-(l-1)$ and $v^{-1}(-l)=w(-(l-1))=-(k-1)$. In this case, we put $\phi(-(k-1),-k)=(-(k-1),-l)$ and $\phi(-(l-1),-l)=(-(l-1),-k)$. Finally, if $k,l\in Z_6$, $k<l$ and $w(k-1)=-(l-1)$, then $v^{-1}(k)=w(k-1)=-(l-1)$ and $v^{-1}(l)=w(l-1)=-(k-1)$. Here we put $\phi(-(k-1),-k)=(-(k-1),l)$ and $\phi(-(l-1),-l)=(-(l-1),k)$.

Thus, we constructed a map $\phi\colon\Au\to\Bu$. One can easily check that if $2\leq-p<-q\leq j$, then
\begin{equation*}
|\Au_i\cap N_{p,q}|\leq|\phi(\Au_i)\cap N_{p,q}|\text{ for all }1\leq i\leq 6.
\end{equation*}
This implies that $\phi$ satisfies (\ref{formula:Au_N}). The proof of case (i) is complete.

ii) Suppose $\beta=\epsi_1+\epsi_j$. This case can be considered similarly to the previous one. (In fact, it is enough to prove that $v^{-1}\geq g_0=us_1$.)

iii) Suppose $\Phi=B_n$ and $\beta=\epsi_1$, or $\Phi=C_n$ and $\beta=2\epsi_1$. In fact, $u=s_{\beta}$. We must prove that $\beta$~divides $d_w$. We will proceed by induction on $n$ (the base $n=2$ is trivial). Note that if $\Phi=B_n$ (resp. $\Phi=C_n$), then $\wt\Phi=B_{n-1}$ (resp. $\wt\Phi=C_{n-1}$). By (\ref{formula_x_v_x_w}a), $l(w)=l(u)+l(v)$ implies $x_w=x_u\cdot x_v$, hence $$c_w=c_{w,\id}=\sum_{g\in\Uu}c_{u,g}\cdot g(c_{v,g^{-1}}),$$ where $\Uu=\{g\in W\mid g\leq u,~g^{-1}\leq v\}$.

Since $u(\alpha_1)=-\epsi_1+\epsi_2<0$, we have $l(us_1)=l(u)-1$. Using (\ref{formula_x_v_x_w}b), we obtain $$c_{u,g}=-\dfrac{c_{us_1,g}+c_{us_1,gs_1}}{g(\alpha_1)}=-\dfrac{c_{us_1,g}}{g(\alpha_1)}-\dfrac{c_{us_1,gs_1}}{g(\alpha_1)}$$ for all $g\in\Uu$. Note that $g(\alpha_1)\neq\beta$, because the length of $\alpha_1=\epsi_1-\epsi_2$ as an element of the Euclidean space $\Rp^n$ is not equal to the length of $\beta$.

Pick an element $g\in\Uu$. Formula (\ref{formula:dyer}) says that there exists a polynomial $h\in S=\Cp[\htt]$ such that $$c_{us_1,g}=h\cdot\prod_{\alpha>0,~s_{\alpha}g\leq u s_1}\alpha^{-1}.$$ Obviously, $us_1=\begin{pmatrix}1&2&3&\ldots&n\\2&-1&3&\ldots&n\end{pmatrix}$. But $g^{-1}\leq v$ implies $g(1)=1$, hence $s_{\beta}g(1)=-1$. By (\ref{formula:Bruhat}), $s_{\beta}g\nleq us_1$. It follows that $\dfrac{c_{us_1,g}}{g(\alpha_1)}=\dfrac{A}{B}$, where $A,B\in S$ are coprime polynomials and $\beta$ does not divide $B$ in $S$.

On the other hand, Formula (\ref{formula:dyer}) claims that there exists a polynomial $f\in S$ such that $$c_{us_1,gs_1}=f\cdot\prod_{\alpha>0,~s_{\alpha}gs_1\leq u s_1}\alpha^{-1}.$$ Recall that $g(1)=1$, because $g^{-1}\leq v$. Assume $g\neq\id$, $s_{\beta}gs_1\leq us_1$ and $g(2)=r$ for some $r$ from $-n$ to $n$, $r\neq1$. Then
\begin{equation*}
\begin{split}
&s_{\beta}gs_1(1)=s_{\beta}g(2)=s_{\beta}(r)=r,\\
&s_{\beta}gs_1(2)=s_{\beta}g(1)=s_{\beta}(1)=-1.
\end{split}
\end{equation*}
Since $us_1(1)=2$ and $s_{\beta}gs_1\leq us_1$, we have $r=2$. Hence if $j\in\{1,2,-2,1\}$, then the $j$th column of $X_{us_1}$ coincides with the $j$th column of $X_{s_{\beta}gs_1}$. Using (\ref{formula:Bruhat}) and proceeding by induction of $j$, one can easily deduce that the $j$th column of the first matrix coincides with the $j$th column of the second one for all $j$. It follows that $X_{us_1}=X_{s_{\beta}gs_1}$, so $s_{\beta}gs_1=us_1$. Hence $g=s_{\beta}u=\id$, a contradiction. Thus, if $g\neq\id$, then $s_{\beta}gs_1\nleq us_1$. This means that if $g\neq\id$, then $\dfrac{c_{us_1,gs_1}}{g(\alpha_1)}=\dfrac{C}{D}$, where $C,D\in S$ are coprime and $\beta$ does not divide $D$.

Now, $g(c_{v,g^{-1}})\in\Cp[\alpha_2,\ldots,\alpha_n]$ for all $g\in\Uu$, so $g(c_{v,g^{-1}})=E/F$, where $E,F\in S$ are coprime and $\beta$ does not divide $F$. We see that $$c_{w}=c_{u,\id}\cdot c_{v,\id}+\dfrac{P}{Q},$$ where $P,Q\in S$ are coprime and $\beta$ does not divide $Q$.
At the same time, $$c_{u,\id}=-\dfrac{c_{us_1,\id}+c_{us_1,s_1}}{\alpha_1}=-\dfrac{c_{us_1,\id}}{\alpha_1}-\dfrac{c_{us_1,s_1}}{\alpha_1}.$$

Since $(s_1us_1)^{-1}(\alpha_1)=s_1u^{-1}s_1(\epsi_1-\epsi_2)=\epsi_1+\epsi_2>0$, we conclude that $l(s_1us_1)=l(us_1)-1$. Hence, by (\ref{formula_x_v_x_w}c),
$$c_{us_1,\id}=\dfrac{s_1(c_{s_1us_1,s_1})-c_{s_1us_1,\id}}{\alpha_1}=-\dfrac{c_{s_1us_1,\id}}{\alpha_1}.$$ The last equality holds, because $s_1us_1(1)=1$, so $s_1us_1\in\wt W$, but $s_1\notin\wt W$, hence $s_1\nleq s_1us_1$ and $c_{s_1us_1,s_1}=0$. Clearly, $c_{s_1us_1,\id}\in\Cp(\alpha_2,\ldots,\alpha_n)$. It turns out that $$-\dfrac{c_{us_1,\id}}{\alpha_1}=\dfrac{c_{s_1us_1,\id}}{\alpha_1^2}=\dfrac{G}{H},$$ where $G,H\in S$ are coprime and $\beta$ does not divide $H$.

Finally, $$c_{us_1,s_1}=\dfrac{s_1(c_{s_1us_1,id})-c_{s_1us_1,s_1}}{\alpha_1}=\dfrac{s_1(c_{s_1us_1,\id})}{\alpha_1}.$$ Note that $s_1us_1=s_{2\epsi_2}$ is an involution in $\wt W$ and $s_1us_1(2)=-2$. By the inductive assumption, if $\Phi=B_n$ (resp. $\Phi=C_n$), then $\wt\beta$ does not divide $\wt d_{s_1us_1}$, where $\wt\beta=\epsi_2$ (resp. $\wt\beta=2\epsi_2$). In other words, $$c_{s_1us_1,\id}=\dfrac{\wt I}{\wt\beta\cdot\wt J},$$ where $\wt I,\wt J\in\wt S=\Cp[\alpha_2,\ldots,\alpha_n]$ are coprime and $\wt\beta$ does not divide $\wt I$. Denote $I=s_1(\wt I)$, $J=s_1(\wt J)$. Since $s_1(\wt\beta)=\beta$, we have $$c_{us_1,s_1}=\dfrac{I}{\alpha_1\cdot\beta\cdot J},$$ where $I,J\in S$ are coprime and $\beta$ does not divide $I$.

Hence there exist coprime polynomials $K,L\in S$ depending only on $u$ such that $$c_{u,\id}=\dfrac{K}{\beta\cdot L}$$ and $\beta$ does not divide $K$. Thus, if $c_{v,\id}=M/N$, $M,N\in\wt S$, then $$c_w=\dfrac{K}{\beta\cdot L}\cdot\dfrac{M}{N}+\dfrac{P}{Q}=\dfrac{KMQ+\beta LNP}{\beta LNQ},$$ where $KMQ\neq0$ and $\beta$ does not divide $KMQ$. Since $d_w=\pm c_w\cdot\prod_{\alpha>0}\alpha$, we conclude that $\beta$ does not divide $d_w$. The proof is complete.}

\sst Things now are ready for the proof of our first main result, Theorem~\ref{mtheo:non_red}. The proof immediately follows from Propositions \ref{prop:non_red_non_eq_C_1} and \ref{prop:non_red_eq_C_1} below (cf. \cite[Propositions 2.6, 2.7, 2.8]{EliseevIgnatyev}). Our goal is to check that if $\sigma,\tau$ are distinct involutions in $W$, then their Kostant--Kumar polynomials do not coincide, and, consequently, the tangent cones $C_{\sigma}$ and $C_{\tau}$ do not coincide as subschemes of $T_p\Fo$. We will proceed by induction on $n$ (the base $n=2$ is trivial).

\propp{Let $\sigma,\tau$ in $W$ be involutions. \label{prop:non_red_non_eq_C_1}If $\Supp{\sigma}\cap\Co_1\neq\Supp{\tau}\cap\Co_1$\textup{,} then $d_{\sigma}\neq d_{\tau}$. In~particular\textup{,} $C_{\sigma}\neq C_{\tau}$.}{We must show that $c_{\sigma}\neq c_{\tau}$. If $\Supp{\sigma}\cap\Co_1=\varnothing$ and $\Supp{\tau}\cap\Co_1=\{\beta\}\neq\varnothing$, then, by Lemma~\ref{lemm:crucial_non_red}, $\beta$ divides $d_{\sigma}$, but $\beta$ does not divide $d_{\tau}$. Hence, $d_{\sigma}\neq d_{\tau}$.

On the other hand, suppose $\Supp{\sigma}\cap\Co_1=\{\beta\}$, $\Supp{\tau}\cap\Co_1=\{\gamma\}$, $\beta\neq\gamma$. We can assume without loss of generality that $s_{\beta}\nleq\tau$. By Lemma~\ref{lemm:crucial_non_red}, $\beta$ does not divide $d_{\sigma}$ and $\gamma$ does not divide~$d_{\tau}$. At the contrary, Formula (\ref{formula:dyer}) shows that there exists a polynomial $f\in S=\Cp[\htt]$ such that $$d_{\tau}=\pm\prod_{\alpha>0}\alpha\cdot c_{\tau}=\pm\prod_{\alpha>0}\alpha\cdot f\cdot\prod_{\alpha>0,~s_{\alpha}\leq\tau}\alpha^{-1}=\pm f\cdot\prod_{\alpha>0,~s_{\alpha}\nleq\tau}\alpha.$$ In particular, $\beta$ divides $d_{\tau}$, thus $d_{\sigma}\neq d_{\tau}$. This completes the proof.}

\propp{Let $\sigma,\tau$ in $W$ be distinct involutions. \label{prop:non_red_eq_C_1}If $\Supp{\sigma}\cap\Co_1=\Supp{\tau}\cap\Co_1$\textup{,} then $d_{\sigma}\neq d_{\tau}$. In~particular\textup{,} $C_{\sigma}\neq C_{\tau}$.}{If $\Supp{\sigma}\cap\Co_1=\Supp{\tau}\cap\Co_1=\varnothing$, then the inductive assumption completes the proof. Suppose $\Supp{\sigma}\cap\Co_1=\Supp{\tau}\cap\Co_1=\{\beta\}$. Let $u$ be as in the proof of Lemma~\ref{lemm:crucial_non_red}. There are three cases:
\begin{equation*}
\begin{split}
&\text{i) }\beta=\epsi_1-\epsi_j,\text{ i.e., }w(1)=j,\\
&\text{ii) }\beta=\epsi_1+\epsi_j,\text{ i.e., }w(1)=-j,\\
&\text{iii) }\Phi=B_n\text{ and }\beta=\epsi_1,\text{ or }\Phi=C_n\text{ and }\beta=2\epsi_1,\text{ i.e., }w(1)=-1.
\end{split}
\end{equation*}
Let us consider these three cases separately.

i) Suppose $\beta=\epsi_1-\epsi_j$. Let $w\in W$ be an involution such that $\Supp{w}\cap\Co_1=\{\beta\}$. Here $u=s_{j-1}\ldots s_2s_1$. Put $g_0=us_1$ and $w=uv$, so $v=u^{-1}w\in\wt W$. Suppose $j>2$. Then we denote $w'=s_{j-1}ws_{j-1}$, $w'=u'v'$ and $h_0=u's_1=s_{j-1}g_0$, where $u'=s_{j-1}u$ and $v'=vs_{j-1}\in\wt W$.

To perform the induction step, we will compare $c_{v,g_0^{-1}}$ with $c_{v',h_0^{-1}}$. It is easy to check that ${v(\alpha_{j-1})<0}$, so $l(vs_{j-1})=l(v)-1$ and, by (\ref{formula_x_v_x_w}b),
\begin{equation*}
c_{v,g_0^{-1}}=-\dfrac{c_{vs_{j-1},g_0^{-1}}+c_{vs_{j-1},g_0^{-1}s_{j-1}}}{g_0^{-1}(\alpha_{j-1})}.
\end{equation*}
Note that $g_0(1)=1$ and $g_0(2)=j$, hence $(R_{g_0})_{j,2}=1$. At the same time, $$(vs_{j-1})^{-1}=s_{j-1}wu=\begin{pmatrix}1&2&\ldots\\1&j-1&\ldots\end{pmatrix},$$
so $(R_{(vs_{j-1})^{-1}})_{j,2}=0$. It follows that $(vs_{j-1})^{-1}\ngeq g_0$, so $vs_{j-1}\ngeq g_0^{-1}$. Thus, $c_{vs_{j-1},g_0^{-1}}=0$ and $$c_{v,g_0^{-1}}=-\dfrac{c_{vs_{j-1},g_0^{-1}s_{j-1}}}{g_0^{-1}(\alpha_{j-1})}=\dfrac{c_{v',h_0^{-1}}}{\epsi_2-\epsi_j}.$$

If $j-1>2$, then we repeat this procedure with $w'$ in place of $w$, etc. In a finite number of steps, we will obtain $w=aw_1a^{-1}$, where $a=s_2s_3\ldots s_{j-1}$, $w_1$ is an involution and $\Supp{w_1}\cap\Co_1=\{\alpha_1\}$. One has $w_1=u_1v_1$, where $u_1=s_1$ and $v_1\in\wt W$ is an involution. Furthermore, $c_{v,g_0^{-1}}=f\cdot c_{v_1,\id}$, where $$f=\dfrac{1}{(\epsi_2-\epsi_j)\cdot(\epsi_2-\epsi_{j-1})\cdot\ldots\cdot(\epsi_2-\epsi_3)}$$
depends only on $j$.

\newpage Now, arguing as in the last two paragraphs of the proof of \cite[Proposition 2.8]{EliseevIgnatyev}, one can conclude the proof. Namely, suppose $\sigma=uv_{\sigma}$ and $\tau=uv_{\tau}$, where $v_{\sigma}$, $v_{\tau}\in\wt W$. Denote $\sigma_1=a\sigma a^{-1}$, $\tau_1=a\tau a^{-1}$, $v_{\sigma}^1=s_1\sigma_1$, and $v_{\tau}^1=s_1\tau_1$. Since $\sigma\neq\tau$, one has $v_{\sigma}^1\neq v_{\tau}^1$. By the inductive assumption, $$c_{v_{\sigma}^1}=\wt c_{v_{\sigma}^1}\neq\wt c_{v_{\tau}^1}=c_{v_{\tau}^1},$$ so $c_{v_{\sigma},g_0^{-1}}=f\cdot c_{v_{\sigma}^1}\neq f\cdot c_{v_{\tau}^1}=c_{v_{\tau},g_0^{-1}}$ and $g_0(c_{v_{\sigma},g_0^{-1}})\neq g_0(c_{v_{\tau},g_0^{-1}})$.

Denote $U_{\sigma}=\{g\in W\mid g\leq u,~g^{-1}\leq v_{\sigma}\}$, then
\begin{equation*}
c_{\sigma}=-\dfrac{c_{us_1,g_0}\cdot g_0(c_{v_{\sigma},g_0^{-1}})}{\beta}-\sum_{g\in U_{\sigma},!g\neq g_0}\dfrac{c_{us_1,g}\cdot g(c_{v_{\sigma},g^{-1}})}{g(\alpha_1)}=\dfrac{A}{B}\cdot\dfrac{P_{\sigma}}{\beta\cdot Q_{\sigma}}+\dfrac{C_{\sigma}}{D_{\sigma}}
\end{equation*}
for some polynomials $A$, $B$, $O_{\sigma}$, $Q_{\sigma}$, $C_{\sigma}$, $D_{\sigma}$. Here $-c_{us_1,g_0}=A/B$, $g_0(c_{v_{\sigma},g_0^{-1}})=P_{\sigma}/Q_{\sigma}$. Similarly,
\begin{equation*}
c_{\tau}=\dfrac{A}{B}\cdot\dfrac{P_{\tau}}{\beta\cdot Q_{\tau}}+\dfrac{C_{\tau}}{D_{\tau}},
\end{equation*}
where $g_0(c_{v_{\tau},g_0^{-1}})=P_{\tau}/Q_{\tau}$. Note that $\beta$ divides neither $A$, nor $D_{\sigma}D_{\tau}$. If $d_{\sigma}=d_{\tau}$, then $$\beta BQ_{\sigma}Q_{\tau}(C_{\sigma}D_{\tau})=AD_{\sigma}D_{\tau}(P_{\tau}Q_{\sigma}-P_{\sigma}Q_{\tau}).$$ This implies that $\beta$ divides $P_{\tau}Q_{\sigma}-P_{\sigma}Q_{\tau}$. But the latter polynomial belongs to the subalgebra generated by $\alpha_2,\ldots,\alpha_n$, so $P_{\tau}Q_{\sigma}-P_{\sigma}Q_{\tau}=0$. Hence $g_0(c_{v_{\sigma},g_0^{-1}})=g_0(c_{v_{\tau},g_0^{-1}})$, a contradiction. Thus, $d_{\sigma}\neq d_{\tau}$.

ii) Suppose $\beta=\epsi_1+\epsi_j$. Let $w\in W$ be an involution such that $\Supp{w}\cap\Co_1=\{\beta\}$. Here $u=s_js_{j+1}\ldots s_{n-1}s_ns_{n-1}\ldots s_2s_1$. Put $g_0=us_1$ and $w=uv$, so $v=u^{-1}w\in\wt W$. Suppose $j<n$. Then we denote $w'=s_jws_j$, $w'=u'v'$ and $h_0=u's_1=s_jg_0$, where $u'=s_ju$ and $v'=vs_j\in\wt W$.

Our goal now is to compare $c_{v,g_0^{-1}}$ with $c_{v',h_0^{-1}}$. It is easy to check that $v(\alpha_j)<0$, so $l(vs_j)=l(v)-1$ and, by (\ref{formula_x_v_x_w}b),
\begin{equation*}
c_{v,g_0^{-1}}=-\dfrac{c_{vs_j,g_0^{-1}}+c_{vs_j,g_0^{-1}s_j}}{g_0^{-1}(\alpha_j)}.
\end{equation*}
Note that $g_0(1)=1$ and $g_0(2)=-j$, hence $(R_{g_0})_{-j,2}=1$. At the same time, $$(vs_j)^{-1}=s_jwu=\begin{pmatrix}1&2&\ldots\\1&-(j+1)&\ldots\end{pmatrix},$$
so $(R_{(vs_j)^{-1}})_{-j,2}=0$. It follows that $(vs_j)^{-1}\ngeq g_0$, so $vs_j\ngeq g_0^{-1}$. Thus, $c_{vs_j,g_0^{-1}}=0$ and $$c_{v,g_0^{-1}}=-\dfrac{c_{vs_j,g_0^{-1}s_j}}{g_0^{-1}(\alpha_j)}=\dfrac{c_{v',h_0^{-1}}}{\epsi_2+\epsi_{j+1}}.$$

If $j+1<n$, then we repeat this procedure with $w'$ in place of $w$, etc. In a finite number of steps, we will obtain $w=aw_1a^{-1}$, where $a=s_{n-1}s_{n-2}\ldots s_j$, $w_1$ is an involution and $\Supp{w_1}\cap\Co_1=\{\alpha_n\}$. Now, $w_1=u_1v_1$, where $u_1=s_ns_{n-1}\ldots s_2s_1$ and $v_1\in\wt W$. Furthermore, $c_{v,g_0^{-1}}=f\cdot c_{v_1,g_1^{-1}}$, where $g_1=u_1s_1$ and $$f=\dfrac{1}{(\epsi_2+\epsi_{j+1})\cdot(\epsi_2+\epsi_{j+2})\cdot\ldots\cdot(\epsi_2+\epsi_n)}.$$

Put $w''=s_nw_1s_n$, $w''=u''v''$, where $u''=s_nu_1$, $v''=v_1s_n\in\wt W$, and $h_1=u''s_1=s_ng_1$. Arguing as above, one can show that
\begin{equation*}
c_{v_1,g_1^{-1}}=\dfrac{c_{v'',h_1^{-1}}}{\gamma},\text{ where }\gamma=-g_0^{-1}(\alpha_n)=
\begin{cases}\epsi_2,&\text{if }\Phi=B_n,\\
2\epsi_2,&\text{if }\Phi=C_n.\\
\end{cases}.
\end{equation*}
Note that $w''$ is an involution and $\Supp{w''}\cap\Co_1=\{\epsi_1-\epsi_n\}$. Applying step i), we see that $w''=bw_2b^{-1}$, where $b=s_2s_3\ldots s_{n-1}$, $\Supp{w_2}\cap\Co_1=\{\alpha_1\}$, $w_2=u_2v_2$, $u_2=s_1$, $v_2\in\wt W$ is an involution and $c_{v_1,g_1^{-1}}=f_1\cdot c_{v_2,\id}$, where $$f_1=\dfrac{1}{(\epsi_2-\epsi_n)\cdot(\epsi_2-\epsi_{n-1})\cdot\ldots\cdot(\epsi_2-\epsi_3)}.$$
Finally, we obtain $c_{v,g_0^{-1}}=f_2\cdot c_{v_2,\id}$, where $f_2=f\cdot\gamma^{-1}\cdot f_1$ depends only on $j$. Now, arguing as in the last two paragraphs of the proof of \cite[Proposition 2.8]{EliseevIgnatyev}, one can conclude the proof.

iii) Suppose $\Phi=B_n$ and $\beta=\epsi_1$, or $\Phi=C_n$ and $\beta=2\epsi_1$. In this case, $u=s_{\beta}$. Recall from the proof of Lemma~\ref{lemm:crucial_non_red} that there exist coprime polynomials $K,L\in S$ such that $\beta$ does not divide $K$ and $$c_{u,\id}=\dfrac{K}{\beta L}.$$ Put $v_{\sigma}=u^{-1}\sigma$, $v_{\tau}=u^{-1}\tau$. Arguing as in step (iii) of the proof of Lemma~\ref{lemm:crucial_non_red}, we deduce that
\begin{equation*}
\begin{split}
&c_{\sigma}=c_{\sigma,\id}=c_{u,\id}c_{v_{\sigma},\id}+\dfrac{P_{\sigma}}{Q_{\sigma}},\\
&c_{\tau}=c_{\tau,\id}=c_{u,\id}c_{v_{\tau},\id}+\dfrac{P_{\tau}}{Q_{\tau}}\\
\end{split}
\end{equation*}
for some $P_{\sigma},Q_{\sigma},P_{\tau},Q_{\tau}\in S$ such that $\beta$ divides neither $Q_{\sigma}$ nor $Q_{\tau}$.

Let $c_{v_{\sigma},\id}=M_{\sigma}/N_{\sigma}$, $c_{v_{\tau},\id}=M_{\tau}/N_{\tau}$ for some $M_{\sigma},N_{\sigma},M_{\tau},N_{\tau}\in\wt S$. Assume $c_{\sigma}=c_{\tau}$. Then
\begin{equation*}
\begin{split}
\dfrac{K}{\beta L}\cdot\dfrac{M_{\sigma}}{N_{\sigma}}+\dfrac{P_{\sigma}}{Q_{\sigma}}&=\dfrac{K}{\beta L}\cdot\dfrac{M_{\tau}}{N_{\tau}}+\dfrac{P_{\tau}}{Q_{\tau}},\\
\dfrac{KM_{\sigma}Q_{\sigma}+\beta LN_{\sigma}P_{\sigma}}{\beta LN_{\sigma}Q_{\sigma}}&=\dfrac{KM_{\tau}Q_{\tau}+\beta LN_{\tau}P_{\tau}}{\beta LN_{\tau}Q_{\tau}}.
\end{split}
\end{equation*}
It follows that $$KQ_{\sigma}Q_{\tau}(M_{\sigma}N_{\tau}-M_{\tau}N_{\sigma})=\beta L(N_{\tau}P_{\tau}-N_{\sigma}P_{\sigma})$$ is divisible by $\beta$. Since $K$, $Q_{\sigma}$ and $Q_{\tau}$ are not divisible by $\beta$, we conclude that $\beta$ divides $M_{\sigma}N_{\tau}-M_{\tau}N_{\sigma}$. But the latter polynomial belongs to $\wt S$, so $$M_{\sigma}N_{\tau}-M_{\tau}N_{\sigma}=0.$$

Thus, $c_{v_{\sigma},\id}=c_{v_{\tau},\id}$. But $v_{\sigma}$ and $v_{\tau}$ are distinct involutions in $\wt W$. The inductive hypothesis guarantees that $c_{v_{\sigma},\id}\neq c_{v_{\tau},\id}$. This contradiction shows that $c_{\sigma}\neq c_{\tau}$. The result follows.}

\sect{Types $A_n$ and $C_n$, reduced case}\label{sect:red}

\sst In this\label{sst:red_coadjoint} Section we will prove our second main result, Theorem~\ref{mtheo:red}. Throughout the Section, we will assume that every $\Phi$ is of type $A_n$ or $C_n$. In this Subsection, we describe connections between tangent cones and coadjoint orbits of $U$, the unipotent radical of the Borel subgroup $B$.

Denote by $\gt$, $\bt$, $\nt$ the Lie algebras of $G$, $B$, $U$ respectively, then $T_p\Fo$ is naturally isomorphic to the quotient space $\gt/\bt$. Using the Killing form on $\gt$, one can identify the latter space with the dual space~$\nt^*$. The group $B$ acts on $\Fo$ by conjugation. Since $p$ is $B$-stable, $B$ acts on the tangent space $T_p\Fo\cong\nt^*$. This action is called \emph{coadjoint}. We denote the result of coadjoint action by $b.\lambda$, $b\in B$, $\lambda\in\nt^*$. In 1962, A.A. Kirillov discovered that orbits of this action play an important role in representation theory of $B$ and $U$, see, e.g., \cite{Kirillov1}, \cite{Kirillov2}.

We fix a basis $\{e_{\alpha},~\alpha\in\Phi^+\}$ of $\nt$ consisting of root vectors. Let $\{e_{\alpha}^*,~\alpha\in\Phi^+\}$ be the dual basis of~$\nt^*$. Let $w\in W$ be an involution. For $\Phi=C_n$, the support of an involution was defined in Subsection~\ref{sst:support} (see Definition~\ref{defi:support}). For $\Phi=A_n$, define the support of $w$ by the following rule. As usual, we identify the set $A_{n-1}^+$ of positive roots with $$\{\epsi_i-\epsi_j,~1\leq i<j\leq n\}\subset\Rp^n.$$ Then $\alpha_1=\epsi_1-\epsi_2$, $\ldots$, $\alpha_{n-1}=\epsi_{n-1}-\epsi_n$ are simple roots. The Weyl group of $A_{n-1}$ is isomorphic to~$S_n$, the symmetric group on $n$ letters. An isomorphism is given by $$s_{\epsi_i-\epsi_j}\mapsto(i,j).$$
Now, if $1\leq i<j\leq n$ and $w(i)=j$, then $\epsi_i-\epsi_j\in\Supp{w}$. Put $$f_w=\sum_{\beta\in\Supp{w}}e_{\beta}^*\in\nt^*.$$

\defi{We say that the $B$-orbit $\Omega_w$ of $f_w$ is \emph{associated} with the involution $w$.}

One can easily check that $\Omega_w\subseteq\rtc{w}$. Further, $\rtc{w}$ is $B$-stable (in fact, the tangent cone to an arbitrary Schubert variety is $B$-stable). Orbits associated with involutions were studied by A.N.~{Pa\-nov}~\cite{Panov} and the second author \cite{Ignatyev1}, \cite{Ignatyev2}, \cite{Ignatyev3}, \cite{Ignatyev4}. In particular, it was shown in \cite[Pro\-po\-si\-tion~4.1]{Ignatyev1} and \cite[Theorem 3.1]{Ignatyev2} that
\begin{equation}
\dim\Omega_w=l(w).\label{formula:dim_Omega}
\end{equation}
Since $\dim\rtc{w}=\dim X_w=l(w)$, we conclude that $\overline{\Omega}_w$, the closure of $\Omega_w$, is an irreducible component of $\rtc{w}$ of maximal dimension. (In fact, $\rtc{w}$ is equidimensional.)

Now, assume that $G'$ is a reductive subgroup of $G''$, $T'$ (resp. $T''$) is a maximal torus of $G'$ (resp. of $G''$), $T'=T''\cap G'$, $B'$ (resp.~$B''$) is a Borel subgroup of $G'$ (resp. of $G''$) containing $T'$ (resp.~$T''$), $B'=B''\cap G'$, and $\Phi'$ (resp. $\Phi''$) is the root system of $G'$ (resp. of $G''$) with respect to $T$ (resp. to~$T''$). We denote by $W'$ (resp. by $W''$) the Weyl group of $\Phi'$ (resp. of $\Phi''$). Denote by $\Fo'=G'/B'$, $\Fo''=G''/B''$ the flag varieties. Put $p'=eB'\in\Fo'$, $p''=eB''\in\Fo''$. Let $U'$ (resp. $U''$) be the unipotent radical of $B'$ (resp. of $B''$), $U'=U''\cap B'$. Denote also by $\gt'$,~$\bt'$,~$\nt'$ the Lie algebras of $G'$, $B'$, $U'$ respectively. Define $\gt''$, $\bt''$, $\nt''$ by the similar way. One can consider the dual space $\nt'^*\cong\gt'/\bt'$ as a~subspace of $\nt''^*\cong\gt''/\bt''$. Hence we can consider $T_{p'}\Fo'$ as a subspace of $T_{p''}\Fo''$.

Pick involutions $w_1,w_2\in W'$. Let $C_i'$ be the reduced tangent cone at the point $p'$ to the Schubert subvariety $X_{w_i}'$ of the flag variety $\Fo'$, $i=1,2$. Similarly, let $C_i''$ be the reduced tangent cone at $p''$ to the Schubert subvariety $X_{w_i}''$ of $\Fo''$, $i=1,2$. Denote by $l'$ (resp. by $l''$) the length function on the Weyl group $W'$ (resp. on $W''$). Assume $C_1'=C_2'$. This implies that $$l'(w_1)=l'(w_2).$$ Note that $C_i'\subseteq C_i''$, hence $B''.C_i'\subseteq C_i''$, $i=1,2$. Denote by $\Omega_{w_i}'\subseteq\nt'^*$ the coadjoint $B'$-orbit associated with the involution $w_i$, $i=1,2$; define $\Omega_{w_i}''$ by the similar way. It follows from Formula (\ref{formula:dim_Omega}) that
\begin{equation*}
\begin{split}
l''(w_i)&=\dim C_i''\geq\dim B''.C_i'\geq\dim B''.\Omega_{w_i}'\\
&=\dim\Omega_{w_i}''=l''(w_i),
\end{split}
\end{equation*}
because $\Omega_{w_i}''=B''.\Omega_{w_i}'$. This implies $l''(w_i)=\dim C_i''=\dim B''.C_i'$. But $C_1'=C_2'$, thus $\dim C_1''=\dim C_2''$. We obtain the following result:
\begin{equation}
\text{if $C_1'=C_2'$, then $l''(w_1)=l''(w_2)$.}\label{formula:if_cones_then_ls}
\end{equation}

\sst In this Subsection, \label{sst:red_A_n}we prove Theorem~\ref{mtheo:red} for $A_n$. Let $W''$ be of type $A_{n+1}$. Let
$$A_{n+1}^+=\{\eta_i-\eta_j,~1\leq i<j\leq n+2\},$$
where $\{\eta_i\}_{i=1}^n$ is the standard basis of $\Rp^{n+2}$. Pick numbers $k_1,k_2$ such that $1\leq k_1<k_2\leq n+2$. Put $P=\{k_1,k_2\}$, $Q=\{1,\ldots,n+2\}\setminus P$, and
\begin{equation*}
\begin{split}
\wt W&=\{w\in W''\mid w(i)=i\text{ for all }i\in P\}\cong S_n,\\
\wt W_2&=\{w\in W''\mid w(i)=i\text{ for all }i\in Q\}\cong S_2, \\
W'&=\{w\in W''\mid w(P)=P,~w(Q)=Q\}=\wt W\times\wt W_2.
\end{split}
\end{equation*}
Let $\Phi'$ (resp. $\wt\Phi$) be the root system of $W'$ (resp. of $\wt W$). Clearly, $\Phi'$ (resp. $\wt\Phi$) is of type $A_{n-1}\times A_1$ (resp. of type $A_{n-1}$). Put $G''=\SL_{n+2}(\Cp)$ and denote by $G'$ (resp. by $\wt G$) the subgroup of $G$ corresponding to~$\Phi'$ (resp. to $\wt\Phi$), then $G'\cong\SL_n(\Cp)\times\SL_2(\Cp)$. Put also
\begin{equation*}
\begin{split}
&A=\{1,\ldots,k_1-1\},\\
&B=\{k_1+1,\ldots,k_2-1\},\\
&C=\{k_2+1,\ldots,n+2\}.
\end{split}
\end{equation*}

Now, let $\Phi=A_{n-1}$. We can assume without loss of generality that $G=\SL_n(\Cp)$. We identify $\Phi$ with $\wt\Phi$ by the map $\epsi_k\mapsto\eta_{k'}$, where
\begin{equation*}
k'=\begin{cases}k,&\text{if }k\leq k_1-1,\\
k+1,&\text{if }k_1\leq k\leq k_2-2,\\
k+2,&\text{if }k_2-1\leq k\leq n.
\end{cases}
\end{equation*}
This identifies $G$ (resp. $W$) with $\wt G$ (resp. with $\wt W$). We denote the image in $\wt W$ of an element $w\in W$ under this identification again by $w$. Let $w\in W$ be an involution. Put $w'=ws_{\eta_{k_1}-\eta_{k_2}}$, then, evidently, $w'$~is an involution in $W'$, and $l'(w')=l(w)+1$.

\lemmp{The length of $w'$ in the \label{lemm:l_w_dash_A_n}Weyl group $W''$ equals $$l''(w')=2(k_2-k_1-1)+4|w(A)\cap C|+l(w)+1.$$}{Clearly, $$l''(w')=l'(w')+\#\{\alpha\in\Phi''^+\setminus\Phi'^+\mid w'(\alpha)<0\}.$$ But
\begin{equation*}
\begin{split}
&\Phi''^+\setminus\Phi'^+=\wt A\cup\wt B\cup\wt C,\text{ where}\\
&\wt A=\{\eta_a-\eta_{k_1},~\eta_a-\eta_{k_2},~a\in A\},\\
&\wt B=\{\eta_{k_1}-\eta_b,~\eta_b-\eta_{k_2},~b\in B\},\\
&\wt C=\{\eta_{k_1}-\eta_c,~\eta_{k_2}-\eta_c,~c\in C\}.
\end{split}
\end{equation*}
For example, $w'(\eta_a-\eta_{k_1})=\eta_{w(a)}-\eta_{k_2}<0$ if and only if $w(a)>k_2$, i.e., $w(a)\in C$. On the other hand, $w'(\eta_a-\eta_{k_2})=\eta_{w(a)}-\eta_{k_1}<0$ if and only if $w(a)>k_1$, i.e., $w(a)\in B$ or $w(a)\in C$. Here we consider $w$ as an element of $\wt W\subset W''$, and, at the same time, as an element of $S_{n+2}$. Hence
\begin{equation*}
\begin{split}
\#\{\alpha\in\wt A\mid w'(\alpha)<0\}&=|w(A)\cap C|+(|w(A)\cap B|+|w(A)\cap C|)\\
&=|w(A)\cap B|+2|w(A)\cap C|.
\end{split}
\end{equation*}

Considering two other cases similarly, one can easily check that
\begin{equation*}
\begin{split}
\#\{\alpha\in\wt B\mid w'(\alpha)<0\}&=(|w(B)\cap A|+|w(B)\cap B|)+(|w(B)\cap B|+|w(B)\cap C|)\\&=|w(B)|+|w(B)\cap B|=|B|+|w(B)\cap B|,\\
\#\{\alpha\in\wt C\mid w'(\alpha)<0\}&=(|w(C)\cap A|+|w(C)\cap B|)+|w(C)\cap A|)\\&=2|w(C)\cap A|+|w(C)\cap B|.\\
\end{split}
\end{equation*}
Since $w$ is an involution, $|w(X)\cap Y|=|X\cap w(Y)|$ for arbitrary subsets $X$, $Y$. Thus,
\begin{equation*}
\begin{split}
\#\{\alpha\in\Phi''^+\setminus\Phi'^+\mid w'(\alpha)<0\}&=|w(A)\cap B|+2|w(A)\cap C|+|B|\\
&+|w(B)\cap B|+2|w(C)\cap A|+|w(C)\cap B|\\
&=|A\cap w(B)|+2|w(A)\cap C|+|B|\\
&+|B\cap w(B)|+2|C\cap w(A)|+|C\cap w(B)|\\
&=|B|+|w(B)|+4|w(A)\cap C|=2|B|+4|w(A)\cap C|\\
&=2(k_2-k_1-1)+4|w(A)\cap C|,
\end{split}
\end{equation*}
because $|w(B)|=|B|=k_2-k_1-1$. The result follows.}

\propp{Let $w_1$, $w_2$ be involutions\label{prop:red_A_n} in the Weyl group $W$ of type $A_{n-1}$, $n\geq2$. If $w_1\neq w_2$, then $\rtc{w_1}\neq\rtc{w_2}$ as subvarieties in $T_p\Fo$.}{Assume $\rtc{w_1}=\rtc{w_2}$. In particular, $$l(w_1)=\dim\rtc{w_1}=\dim\rtc{w_2}=l(w_2).$$ Since $w_1\neq w_2$, there exists $1\leq k\leq n$ such that $w_1(\epsi_i)=w_2(\epsi_i)$ for $1\leq i\leq k-1$, and $$w_1(\epsi_k)=\epsi_{m_1}\neq\epsi_{m_2}=w_2(\epsi_k).$$

Assume without loss of generality that $m_1>m_2$. Note that $m_2\geq k$, so $m_1>k$. Let $G'$, $G''$ etc. be as above, where $$1<k_1=k+1<k_2=m_1+1<n+2.$$ Then
$$w_1(A)\cap C=(w_2(A)\cap C)\sqcup\{k\},$$ so, by the previous Lemma, $l''(w_1')\neq l''(w_2')$. On the other hand, $$C_i'=\rtc{w_i}\times\Cp e_{\eta_{k_2}-\eta_{k_1}},$$ $i=1,2$, so $C_1'=C_2'$. This contradicts (\ref{formula:if_cones_then_ls}). The result follows.}

\sst In this Subsection, \label{sst:red_C_n}we prove Theorem~\ref{mtheo:red} for $C_n$. As in Section~\ref{sect:non_red}, we identify $C_n^+$ with $$\{\epsi_i-\epsi_j,~\epsi_i+\epsi_j,~1\leq i<j\leq n\}\cup\{2\epsi_i,~1\leq i\leq n\}\subset\Rp^n$$
and the Weyl group $W$ of type $C_n$ with the subgroup of $S_{\pm n}$ consisting of $w$ such that $w(-i)=-w(i)$ for all $1\leq i\leq n$.

Let $W''$ be of type $C_{n+2}$. We identify $C_{n+2}^+$ with
$$\{\eta_i-\eta_j,~\eta_i+\eta_j,~1\leq i<j\leq n+2\}\cup\{2\eta_i,~1\leq i\leq n+2\}.$$
Pick numbers $k_1,k_2$ such that $1\leq k_1<k_2\leq n+2$. As in the previous Subsection, put $P=\{k_1,k_2\}$, $Q=\{1,\ldots,n+2\}\setminus P$, and
\begin{equation*}
\begin{split}
\wt W&=\{w\in W''\mid w(i)=i\text{ for all }i\in P\},\\
\wt W_2&=\{w\in W''\mid w(i)=i\text{ for all }i\in Q\}, \\
W'&=\{w\in W''\mid w(P)=P,~w(Q)=Q\}.
\end{split}
\end{equation*}
Let $\Phi'$ (resp. $\wt\Phi$) be the root system of $W'$ (resp. of $\wt W$). Clearly, $\Phi'$ (resp. $\wt\Phi$) is of type $C_n\times C_2$ (resp.~$C_n$). Put $G''=\Sp_{2n+4}(\Cp)$ and denote by $G'$ (resp. by $\wt G$) the subgroup of $G$ corresponding to~$\Phi'$ (resp. to $\wt\Phi$), then $G'\cong\Sp_{2n}(\Cp)\times\Sp_4(\Cp)$. As above, put also
\begin{equation*}
A=\{1,\ldots,k_1-1\},~B=\{k_1+1,\ldots,k_2-1\},~C=\{k_2+1,\ldots,n+2\}.
\end{equation*}

Now, let $\Phi=C_n$. We can assume without loss of generality that $G=\Sp_{2n}(\Cp)$. We identify $\Phi$ with~$\wt\Phi$ by the map $\epsi_k\mapsto\eta_{k'}$, where
\begin{equation*}
k'=\begin{cases}k,&\text{if }k\leq k_1-1,\\
k+1,&\text{if }k_1\leq k\leq k_2-2,\\
k+2,&\text{if }k_2-1\leq k\leq n.
\end{cases}
\end{equation*}
This identifies $G$ (resp. $W$) with $\wt G$ (resp. with $\wt W$). We denote the image in $\wt W$ of an element $w\in W$ under this identification again by $w$. For any $X\subseteq\{1,\ldots,n+2\}$, put $X^-=-X$ and $X^{\pm}=X\cup X^-$. Let $w\in W$ be an involution. Arguing as in the proof of Lemma~\ref{lemm:l_w_dash_A_n}, we obtain the following result.

\mlemm{\textup{i)} If\label{lemm:l_w_dash_C_n} $w'=ws_{\eta_{k_1}-\eta_{k_2}-1}$\textup{,} then $$l''(w')=2(k_2-k_1)+4|w(A)\cap C^{\pm}|+4|w(A)\cap B^-|+4|w(A)\cap A^-|+l(w)+1.$$
\textup{ii)} If $w'=ws_{\eta_{k_1}+\eta_{k_2}}$\textup{,} then $$l''(w')=2(k_2-k_1-1)+4|C|+4|w(A)\cap B^-|+4|w(A)\cap A^-|+l(w)+3.$$}

\propp{Let $w_1$, $w_2$ be involutions\label{prop:red_C_n} in the Weyl group $W$ of type $C_n$. If $w_1\neq w_2$, then $\rtc{w_1}\neq\rtc{w_2}$ as subvarieties in $T_p\Fo$.}{Assume $\rtc{w_1}=\rtc{w_2}$. In particular, $$l(w_1)=\dim\rtc{w_1}=\dim\rtc{w_2}=l(w_2).$$ Since $w_1\neq w_2$, there exists $1\leq k\leq n$ such that $w_1(\epsi_i)=w_2(\epsi_i)$ for $1\leq i\leq k-1$, and $w_1(\epsi_k)\neq w_2(\epsi_k)$ (signs are independent).

We can assume without loss of generality that $w_1(\epsi_k)<w_2(\epsi_k)$, i.e, $w_2(\epsi_k)-w_1(\epsi_k)$ is a sum of positive roots. Note that $w_1(\epsi_k)\neq\epsi_k$. Put $k_1=k+1$ and consider four different cases.

i) Suppose $w_1(\epsi_k)=-\epsi_k$. Then $w_2(\epsi_k)=\pm\epsi_l$ for some $l>k$. Put $k_2=k_1+1$, so $B=\varnothing$ and $$w_1(A)\cap A^-=(w_2(A)\cap A^-)\sqcup\{k\},$$ hence, by Lemma~\ref{lemm:l_w_dash_C_n} (ii), $l''(w_1')\neq l''(w_2')$, where $w_i'=w_is_{\eta_{k_1}+\eta_{k_2}}$, $i=1,2$. On the other hand, $\rtc{w_1}=\rtc{w_2}$ implies $C_1'=C_2'$. This contradicts (\ref{formula:if_cones_then_ls}).

ii) Next, suppose $w_1(\epsi_k)\neq-\epsi_k$, $w_1(\epsi_k)<0$, $w_2(\epsi_k)>0$. Put $k_2=n+2$, so $C=\varnothing$ and $$(w_i(A)\cap B^-)\cup(w_i(A)\cap A^-)=w_i(A)\cap\{-1,\ldots,-n\},~i=1,2.$$ Hence, by Lemma~\ref{lemm:l_w_dash_C_n} (i), $l''(w_1')\neq l''(w_2')$, where $w_i'=w_is_{\eta_{k_1}-\eta_{k_2}}$, $i=1,2$. On the other hand, $C_1'=C_2'$. This contradicts (\ref{formula:if_cones_then_ls}).

\newpage iii) Now, suppose $w_1(\epsi_k)=\epsi_{m_1}$, $w_2(\epsi_k)=\epsi_{m_2}$, $m_1>m_2$. Put $k_2=m_1+1$, then $w_1(k)\in C$ and $w_2(k)\in B$ (here we consider $w_1$ and $w_2$ as an elements of $\wt W\subseteq W''$, or, equivalently, as an elements of~$S_{\pm n}$). Lemma~\ref{lemm:l_w_dash_C_n} (i) shows that $l''(w_1')\neq l''(w_2')$, where $w_i'=w_is_{\eta_{k_1}-\eta_{k_2}}$, $i=1,2$, but $C_1'=C_2'$, a~contradiction.

iv) Finally, suppose $w_1(\epsi_k)\neq-\epsi_k$, $w_1(\epsi_k)=-\epsi_{m_1}$, $w_2(\epsi_k)=-\epsi_{m_2}$, $m_1>m_2$. As above, put $k_2=m_1+1$, then $w_1(k)\in C^-$ and $w_2(k)\in B^-$. Lemma~\ref{lemm:l_w_dash_C_n} (ii) says that $l''(w_1')\neq l''(w_2')$, where $w_i'=w_is_{\eta_{k_1}+\eta_{k_2}}$, $i=1,2$, but $C_1'=C_2'$, a contradiction. This completes the proof.}

\end{document}